\documentclass[reqno]{amsart}
\usepackage{mathrsfs}
\usepackage{color}
\usepackage{graphics}
\usepackage{v2divide}

\def\minimax#1#2#3{%
	\def\pnamem{minimum}
	\def\pnameM{maximum}
	\def\anameU{_up}
	\def\anameD{_down}
	\def\anameS{}
	\ifx#2m
		\raise 6truept\hbox{$\scriptstyle #1$}%
		\kern -2.5truept
		\ifx#3U\else\kern -1truept\fi
	\fi
	\vcenter{\hbox{\noindent
		\input{x-\expandafter\csname pname#2\endcsname%
		\csname aname#3\endcsname.pstex_t}}}%
	\kern -3truept
	\ifx#2M
		\kern -3.5truept
		\ifx#3U\else\kern -0.5truept\fi
		\raise 6truept\hbox{$\scriptstyle #1$}%
	\fi
}

\def\vertexup#1{%
	\kern 2.5truept
	\lower 6truept\hbox{$\scriptstyle #1$}%
	\kern -8truept
	\vcenter{\hbox{\noindent \begin{picture}(0,0)%
\includegraphics{vertex_up.pstex}%
\end{picture}%
\setlength{\unitlength}{3947sp}%
\begingroup\makeatletter\ifx\SetFigFont\undefined
\def\x#1#2#3#4#5#6#7\relax{\def\x{#1#2#3#4#5#6}}%
\expandafter\x\fmtname xxxxxx\relax \def\y{splain}%
\ifx\x\y   
\gdef\SetFigFont#1#2#3{%
  \ifnum #1<17\tiny\else \ifnum #1<20\small\else
  \ifnum #1<24\normalsize\else \ifnum #1<29\large\else
  \ifnum #1<34\Large\else \ifnum #1<41\LARGE\else
     \huge\fi\fi\fi\fi\fi\fi
  \csname #3\endcsname}%
\else
\gdef\SetFigFont#1#2#3{\begingroup
  \count@#1\relax \ifnum 25<\count@\count@25\fi
  \def\x{\endgroup\@setsize\SetFigFont{#2pt}}%
  \expandafter\x
    \csname \romannumeral\the\count@ pt\expandafter\endcsname
    \csname @\romannumeral\the\count@ pt\endcsname
  \csname #3\endcsname}%
\fi
\fi\endgroup
\begin{picture}(174,174)(589,77)
\end{picture}
}}%
	\kern -3.5truept
}

\def\vertexdown#1{%
	\kern 3truept
	\raise 6truept\hbox{$\scriptstyle #1$}%
	\kern -8.5truept
	\vcenter{\hbox{\noindent \begin{picture}(0,0)%
\includegraphics{vertex_down.pstex}%
\end{picture}%
\setlength{\unitlength}{3947sp}%
\begingroup\makeatletter\ifx\SetFigFont\undefined
\def\x#1#2#3#4#5#6#7\relax{\def\x{#1#2#3#4#5#6}}%
\expandafter\x\fmtname xxxxxx\relax \def\y{splain}%
\ifx\x\y   
\gdef\SetFigFont#1#2#3{%
  \ifnum #1<17\tiny\else \ifnum #1<20\small\else
  \ifnum #1<24\normalsize\else \ifnum #1<29\large\else
  \ifnum #1<34\Large\else \ifnum #1<41\LARGE\else
     \huge\fi\fi\fi\fi\fi\fi
  \csname #3\endcsname}%
\else
\gdef\SetFigFont#1#2#3{\begingroup
  \count@#1\relax \ifnum 25<\count@\count@25\fi
  \def\x{\endgroup\@setsize\SetFigFont{#2pt}}%
  \expandafter\x
    \csname \romannumeral\the\count@ pt\expandafter\endcsname
    \csname @\romannumeral\the\count@ pt\endcsname
  \csname #3\endcsname}%
\fi
\fi\endgroup
\begin{picture}(174,174)(589,77)
\end{picture}
}}%
	\kern -3truept
}

\begin{document}
\title[Casson invariant of divide knots]{Casson invariant of knots associated
with divides}
\author[A.~Shumakovitch]{Alexander Shumakovitch}
\address{Mathematisches Institut, Universit\"at Basel, Rheinsprung 21,
CH-4051, Basel,\hfill\break Switzerland}
\email{Shurik@math.unibas.ch}
\thanks{The author is partially supported by the Swiss National Science
Foundation}
\subjclass{57M25, 57M27}
\keywords{Divide, Casson invariant, Arnold's invariants of plane curves}
\begin{abstract}
I present a formula for the Casson invariant of knots associated with divides.
The formula is written in terms of Arnold's invariants of pieces of the
divide. Various corollaries are discussed.
\end{abstract}
\maketitle

\section*{Introduction}
A divide $P$ is the image of a generic relative immersion of a $1$-dimensional
compact (not necessarily connected) manifold into the standard unit disc
$D=\{(x,y)\in\R^2\mid x^2+y^2\le1\}$ in the plane $\R^2$. The relativity of
the immersion implies that the boundary of the manifold in mapped into the
boundary of $D$. The genericity condition means that a divide has only
transversal double points as singularities and is transversal to $\p D$ at its
boundary points. Divides are considered modulo ambient isotopy in $D$.

With every divide $P$ one can associate a link $L(P)\subset S^3$ using the
following construction. Consider a standard projection $\pi$ of
$S^3=\{(x,y,u,v)\in\R^4\mid x^2+y^2+u^2+v^2=1\}$ onto $D$, i.e.
$\pi(x,y,u,v)=(x,y)$. Preimage of a point $p\in D$ under $\pi$ is either a
circle, if $p$ belongs to the interior of $D$, or a single point.

Now $L(P)$ consists of all points $(x,y,u,v)\in\pi^{-1}(P)\subset S^3$ such
that either $(x,y)\in\p P\subset\p D$ or $(u,v)$ is tangent to $P$ at $(x,y)$.
It is easy to see that for any $p\in P$ its preimage $\pi^{-1}(p)\cap L(P)$
in $L(P)$ consists of either $1$, $4$ or $2$ points, depending on whether $p$
is a boundary, double or generic point of $P$. $L(P)$ is indeed a link with
$2c+i$ components, where $c$ and $i$ are the numbers of closed and non-closed
components of $P$, respectively. Ambiently isotopic divides obviously give
rise to ambiently isotopic links.

Divides and associated links were originally considered by N.~A'Campo
\cite{ACampo-divide-def} and are closely related to the real morsifications of
isolated complex plane curve singularities (see also
\cite{ACampo-divide-1,ACampo-divide-2}).

In this paper I present a formula for the Casson invariant of the knot
associated with a given divide with no closed and only one non-closed
components. Such divides are called {\em I-divides} in this paper. The Casson
invariant of a link $L$ can be defined as $\frac12\GD''_L(1)$, where
$\GD_L(t)$ is the Alexander polynomial of $L$ and $\GD''_L(1)$ is the value of
its second derivative at $1$. It is also a unique Vassiliev invariant of
degree $2$ that takes values $0$ on the unknot and $1$ on a trefoil. Since
the Casson invariant of a link $L$ is equal to the sum of its values on the
components of $L$, only the case of knots is interesting.

Amazingly enough, the formula is written in terms of Arnold's invariants
$J^\pm$ and $\St$ of pieces of the divide (see section~\ref{sec:arnold}
or~\cite{Arnold-curve-invars} for definitions).

In \cite{Chmutov-J_pm} S.~Chmutov constructed a second order $J^\pm$-type
invariant of long curves, called $J^\pm_2$, and proved that its value on an
I-divide is equal to the Casson invariant of the associated knot. Since the
invariant was defined by an actuality table only, the computation of its
values on a curve with even as few as $2$ double points was rather involved.
The original definition also did not provide a mean to compute changes of
$J^\pm_2$ under self-tangency perestroikas. My formula fills these gaps (see
Lemmas~\ref{lem:inverse-move},~\ref{lem:direct-move} and
section~\ref{sec:chmutov}).

\subsection*{Acknowledgements} I am deeply grateful to Norbert~A'Campo for his
valuable remarks and encouragement.

\section{Main ingredients and definitions}

\subsection{Diagrams of divide links}
Let $P\subset D$ be a divide. In this section I present a way to draw a
usual diagram with over- and under-crossings of the link $L(P)$. The following
construction is due to M.~Hirasawa~\cite{Hirasawa-diagram} with some minor
modifications (see also \cite{Gibson-Ishikawa-diagram}).

One first perturbs $P$ slightly without changing its ambient isotopy type so
that the following two conditions are met:
\begin{itemize}
\item at every double point the two branches of $P$ are parallel to the main
diagonals $y=\pm x$ of the plane;
\item $P$ has only finitely many points where the tangent vector is parallel
to the $y$-axis, and at every such a point the projection of $P$ onto the
$x$-axis has either a local minimum or a local maximum.
\end{itemize}

\begin{figure}
\centerline{\input{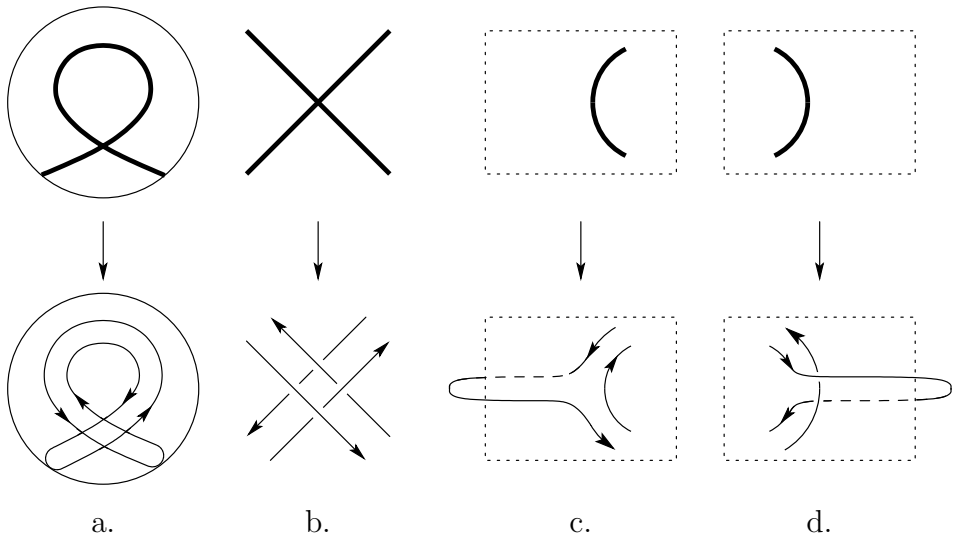}}
\caption{Drawing a diagram of the divide link}
\label{fig:divide-link}
\end{figure}

\begin{figure}
\centerline{\input{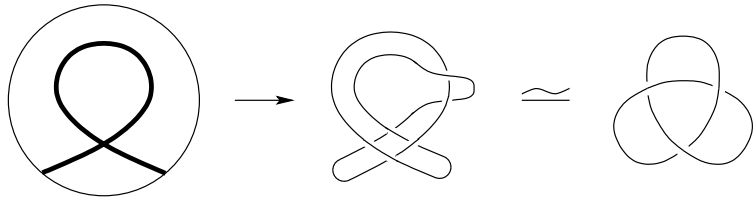}}
\caption{Example of a divide link}
\label{fig:divide-example}
\end{figure}

The next step is to double the divide, i.e. to draw two parallel strands along
the divide and to join them at the end points (see
Figure~\ref{fig:divide-link}.a). The resulting closed curves have a natural
orientation according to the {\em ``keeping right''\/} rule. One puts over-
and under-crossings near double points according to the rule depicted in
Figure~\ref{fig:divide-link}.b. And finally at every point of $P$ with the
tangent vector parallel to the $y$-axis the strand pointing downwards should
make a {\em ``jump through infinity''\/}. That means to go far right (left)
above (below) the rest of the diagram and to return back below (above) it,
according to whether the projection of $P$ onto the $x$-axis has a local
maximum or minimum at that point (see Figure~\ref{fig:divide-link}.c,d).
Denote the diagram obtained by $D(P)$. For example, the divide depicted in
Figure~\ref{fig:divide-link}.a gives rise to a (right) trefoil (see
Figure~\ref{fig:divide-example}).

\begin{thm}[Hirasawa \cite{Hirasawa-diagram}]
The diagram $D(P)$ represents the link $L(P)$ in $S^3$.
\end{thm}

\subsection{Arnold's invariants of plane curves}\label{sec:arnold}
Let $C$ be a generic plane curve, i.e. a ($C^1$-smooth) immersion of the
circle into the plane that has only transversal double points as
singularities. For any such a curve one can define its {\em Whitney index\/}
or simply {\em index\/} as the total rotation number of tangent vector
to the curve. This number is, clearly, the degree of the map that associates
a direction of the tangent vector to every point of the circle. The index of
a curve $C$ is denoted by $w(C)$. It is easy to see that the index does not
change under a {\em regular homotopy\/} of a curve that is a $C^1$-smooth
homotopy in the class of $C^1$-immersions. Moreover, two plane curves are
regular homotopic if and only if their Whitney indices are the same.

\begin{figure}
\centerline{\input{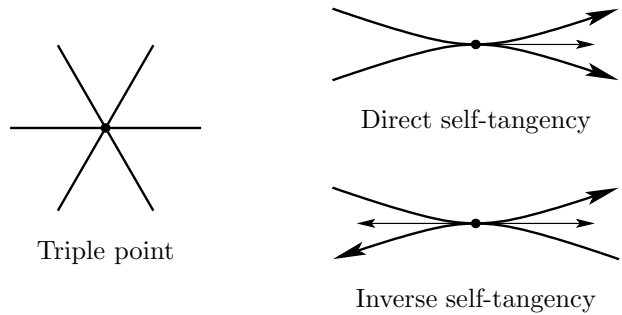}}
\caption{Three types of singularities that may appear in a generic regular
homotopy}
\label{fig:singular-points}
\end{figure}

In a generic regular homotopy connecting two generic curves only a finite
number of non-generic ones can appear. Each of these curves differs from a
generic one either in exactly one point of triple transversal
self-intersection or in exactly one point of self-tangency. In a point of
self-tangency the velocity vectors of the tangent branches can have either the
same directions or the opposite ones. In the first case the self-tangency is
said to be {\em direct\/} and in the second {\em inverse\/} (see
Figure~\ref{fig:singular-points}). The type of the self-tangency does not change
under reversing of orientation.

Hence there are three types of singular curves that may appear in a generic
regular homotopy. Passages through these curves correspond to three {\em
perestroikas\/} of generic curves.

Consider the triple point perestroika more carefully. Just before and just
after the passage through a singular curve with a triple point, there is a
small triangle close to the place of perestroika, which is formed by three
branches of the curve. This triangle is said to be {\em vanishing\/}. The
orientation of the curve defines a cyclic order of sides of the triangle. This
is the order in which one meets the sides while traveling along the curve.
This cyclic order gives the orientation of the triangle and, therefore, the
orientation of its sides. Denote by $q$ the number of sides of the vanishing
triangle for which the orientation obtained coincides with the orientation of
the curve ($q$, obviously, takes value between $0$ and $3$).

Define a {\em sign\/} of the vanishing triangle as $(-1)^q$. The sign does
not change under reversing of orientation of the curve. One can easily check
that before and after the perestroika the sings of the vanishing triangles are
different.

\begin{defins}[Arnold \cite{Arnold-curve-invars}]
1.~A triple point perestroika is said to be {\em positive\/} if
the newborn vanishing triangle is positive.

2.~A self-tangency perestroika is said to be {\em positive\/} if
it increases (by $2$) the number of self-intersection points of the curve.
\end{defins}

\begin{figure}
\centerline{\input{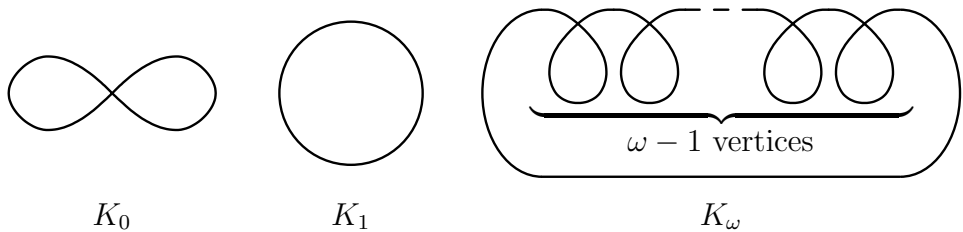}}
\caption{The standard curves with Whitney indices $0,\pm 1,\pm \Go$}
\label{fig:standard-curves}
\end{figure}

The following theorem provides a definition of invariants of generic
(plane) curves.
\begin{thm}[Arnold \cite{Arnold-curve-invars}]
There exist three integers $\St(C)$, $J^+(C)$, and $J^-(C)$ assigned
to an arbitrary generic plane curve $C$ that are uniquely defined by
the following properties.
\begin{enumerate}
\renewcommand{\labelenumi}{(\roman{enumi})}
\item $\St$, $J^+$ and $J^-$ are invariant under a regular homotopy in the
class of generic curves.
\item $\St$ does not change under self-tangency perestroikas and
{\em increases\/} by $1$ under a positive triple point perestroika.
\item $J^+$ does not change under triple point and inverse self-tangency
perestroikas and {\em increases\/} by $2$ under a positive direct
self-tangency perestroika.
\item $J^-$ does not change under triple point and direct self-tangency
perestroikas and {\em decreases\/} by $2$ under a positive inverse
self-tangency perestroika.
\item On the standard curves $K_\Go$, shown in
Figure~\ref{fig:standard-curves}, $\St$, $J^+$, and $J^-$ take the following
values\rm:
\begin{alignat*}{3}
\St(K_0)&=0,  &\qquad \St(K_{\Go+1})&=\Go &&\quad
(\Go=0, 1, 2, \dotsc);\\
J^+(K_0)&=0,  &\qquad J^+(K_{\Go+1})&=-2\Go &&\quad
(\Go=0, 1, 2, \dotsc);\\
J^-(K_0)&=-1, &\qquad J^-(K_{\Go+1})&=-3\Go &&\quad
(\Go=0, 1, 2, \dotsc).
\end{alignat*}
\end{enumerate}
\end{thm}

\begin{rem}
The normalization of $\St$ and $J^\pm$, which is fixed by the last property,
makes them additive with respect to the connected sum of curves. It is easy to
see that the invariants are independent of the orientation of a curve.
\end{rem}

\subsection{Formulas for Arnold's invariants}
It is very complicated to use the definition given above to compute Arnold's
invariants for a given curve with sufficiently many double points. The
two theorems stated below provide explicit formulas, which make those
computations easy and straightforward.

Consider a generic plane curve $C$. For every region $r$ of $C$ one defines
its {\em index\/} with respect to $C$ as the total rotation number of the
radius vector that connects an arbitrary interior point of $r$ to a point
traveling along $C$. This number is clearly independent of the choice of the
point in $r$ and is denoted by $\ind_C(r)$. An {\em index\/} $\ind_C(e)$ of an
edge $e$ of $C$ is defined as a half-sum of the indices of the two regions
adjacent to $e$. An {\em index\/} $\ind_C(v)$ of a double point $v$ of $C$ is a
quarter-sum of the indices of the four regions adjacent to $v$.

\begin{figure}
\centerline{\input{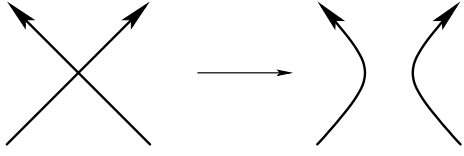}}
\caption{The curve smoothing at a self-intersection point}
\label{fig:curve-smoothing}
\end{figure}

\begin{thm}[Viro \cite{Viro-J's-formulas}]
Let $C$ be a generic plane curve and $\tC$ be a family of embedded
circles obtained as a result of smoothing of the curve $C$ at each double
point with respect to the orientation (see Figure~\ref{fig:curve-smoothing}).
Then
\begin{align}
J^+(C)&=1-\sum_{r\in\ScrR_{\tC}}\ind^{\, 2}_{\tC}(r)\,\chi(r)+n,\\
J^-(C)&=1-\sum_{r\in\ScrR_{\tC}}\ind^{\, 2}_{\tC}(r)\,\chi(r),
\end{align}
where $\ScrR_{\tC}$ is the set of all regions of $\tC$, $\chi$ is the
Euler characteristic, and $n$ is the number of double points of the curve $C$.
\end{thm}

Fix now a base point $f$ on a generic plane curve $C$ that is not a double
point. One can enumerate all edges by numbers from $1$ to $2n$ (where $n$ is
again the number of double points of $C$) following the orientation and
assigning $1$ to the edge with the point $f$.

Consider an arbitrary double point $v$ of $C$. There are two edges pointing to
$v$. Let them have numbers $i$ and $j$ such that the tangent vector to the
edge $i$ and the tangent vector to the edge $j$ give a positive orientation of
the plane. One defines a {\em sign\/} $s(v)$ of $v$ to be $\sign(i-j)$.

\begin{thm}[\cite{Myself-strangeness}]
Let $C$ be a generic plane curve. Then
\begin{equation}
\St(C)=\sum_{v\in\ScrV_C}\ind_C(v)s(v)+\ind^{\,2}_C(f)-\frac14,
\end{equation}
where $\ScrV_C$ is the set of all double points of $C$.
\end{thm}

\begin{rem}
If the point $f$ is chosen on an exterior edge (that is an edge that bounds
the exterior region), then $\ind_C(f)=\pm\frac12$ and one gets
\begin{equation}
\St(C)=\sum_{v\in\ScrV_C}\ind_C(v)s(v).
\end{equation}
\end{rem}

\begin{rem}
One can find another formulas for Arnold's invariants in
\cite{Chmuzhin-invars,Polyak-arnold-gauss,Myself-strangeness}.
\end{rem}

\section{Formula for the Casson invariant}
\subsection{The main results}
Let $P\subset D$ be an I-divide. Choose an arbitrary orientation on $P$. For
any double point $v$ of $P$ denote by $O_v$ and $I_v$ closed and non-closed
curves\footnote{Notation rationale: letters $O$ and $I$ were chosen because
they conveniently represent the topological type of the corresponding curves.}
obtained by smoothing $P$ at $v$ with respect to the orientation (see
Figure~\ref{fig:curve-smoothing}). The end points of $P$ split the boundary
$\p D$ into two arcs. Complete $P$ with one of those arcs, so that the
orientation of $P$ induces a counter-clockwise orientation on the arc. Call
the resulting closed curve $\overline P$. Finally, for a given plane curve
$C$ denote $1-J^-(C)$ by $\tJ(C)$.

\begin{figure}
\centerline{\input{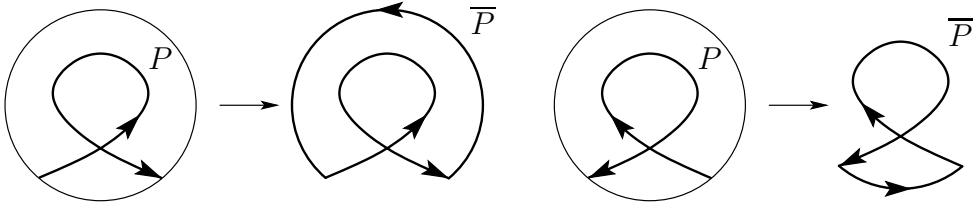}}
\caption{Closure of a divide depending on the orientation}
\label{fig:divide-closure}
\end{figure}

\begin{rem}
The closure $\overline P$ of $P$ clearly depends on the orientation of $P$
(see Figure~\ref{fig:divide-closure}).
\end{rem}

\begin{rem}
$1-J^-(C)$ can be expressed with a very simple formula
$1-J^-(C)=\sum_{r\in\ScrR_{\tC}}\ind^{\, 2}_{\tC}(r)\,\chi(r)$ and therefore
deserves a separate notation.
\end{rem}

\begin{thm}\label{thm:casson-formula}
Let $P$ be an I-divide (i.e. a divide with only one component that is
non-closed). Then the Casson invariant $v_2$ of the knot $L(P)$ is given by
\begin{equation}\label{eq:casson-formula}
v_2(L(P))=\sum_{v\in\ScrV_P}(\tJ(O_v)+\frac14\#(O_v\cap
I_v))+\frac{J^+(\overline P)+2\St(\overline P)}4,
\end{equation}
where $\ScrV_P$ is the set of all double points of $P$ and $\#(O_v\cap
I_v)$ is the number of intersection points of $O_v$ with $I_v$.
\end{thm}

\begin{rem}
$J^+(\overline P)$ and $\St(\overline P)$ obviously depend on the orientation
of $P$. But the sum $J^+(\overline P)+2\St(\overline P)$ does not, since all
the other ingredients of \eqref{eq:casson-formula} are independent of the
orientation. This fact can also be verified directly.
\end{rem}

\subsection{Special cases}
It is well known~\cite{Aicardi-tree-curves} that for a tree-like generic plane
curve $C$ the sum $J^+(C)+2\St(C)$ is always $0$. Here a curve $C$ is called
{\em tree-like\/} if the smoothing at every double point of $C$ produces two
disjoint curves. One can define a {\em tree-like divide\/} in a similar
fashion. It follows from the definition that $\#(O_v\cap I_v)=0$ for every
double point $v$ of a tree-like divide $P$. Hence the
formula~\eqref{eq:casson-formula} can be simplified.

\begin{prop}
Let $P$ be a tree-like divide. Then
\begin{equation}
v_2(L(P))=\sum_{v\in\ScrV_P}\tJ(O_v).
\end{equation}
\end{prop}

A slalom divide (see~\cite{ACampo-divide-3} for definitions) is a special case
of tree-like divides. In this case the numbers $\tJ(O_v)$ can be
calculated directly.

\begin{prop}
Let $P$ be a slalom divide. Then
\begin{equation}
v_2(L(P))=\sum_{v\in\ScrV_P}(1+n(O_v)),
\end{equation}
where $n(O_v)$ is the number of the double points of $O_v$. Since this number
is always non-negative, the Casson invariant of a slalom knot is always
strictly positive.
\end{prop}

\begin{rem}
Formula~\ref{eq:casson-formula} admits an easy rewrite in terms of Gauss
diagrams of divides (see~\cite{Polyak-arnold-gauss} for definitions and
examples of the technique). In particular, the term
$\displaystyle\sum_{v\in\ScrV_P}\frac14\#(O_v\cap I_v)$ is half the number of
generic intersections of chords in the Gauss diagram of a divide $P$.
\end{rem}

\section{Proof of the Theorem~\ref{thm:casson-formula}}

\subsection{Normalization}
Denote the right-hand side of~\eqref{eq:casson-formula} by $\ScrX(P)$. It is
easy to check that the desired equality $v_2(L(P))=\ScrX(P)$ holds true for
the standard divides $D_0$, $D_1,\,\dots$\,. Here a standard divide $D_k$ looks
like a standard curve $K_{k+1}$ cut at the external edge and appropriately
immersed into the unit disc. Indeed, $v_2(L(D_0))=0=\ScrX(D_0)$ and
$v_2(L(D_1))=1=\ScrX(D_1)$, since $L(D_1)$ is a trefoil knot (see
Figure~\ref{fig:divide-example}). Furthermore, $L(D_n)$ is the connected sum
of $n$ copies of $L(D_1)$ for $n>1$ and, hence, $v_2(L(D_n))=v_2(\#nL(D_1))= n
v_2(L(D_1))=n=n\tJ(K_1)=\sum_{v\in\ScrV_{D_n}}\tJ(O_v)=\ScrX(D_n)$.

\subsection{Invariance under perestroikas}
Since one can transform any divide into a standard one with a finite sequence
of Arnold's perestroikas, it is enough to check that $v_2(L(P))$ and
$\ScrX(P)$ change in the same way under the perestroikas. It is proved in
Lemmas which follow.

\begin{lem}\label{lem:inverse-move}
Let $P'$ be the result of a positive inverse self-tangency perestroika of $P$.
Assign $\GD\ScrX=\ScrX(P')-\ScrX(P)$ and $\GD v_2=v_2(L(P'))-v_2(L(P))$. Then
$\GD\ScrX=\GD v_2$.
\end{lem}

\begin{figure}
\centerline{\input{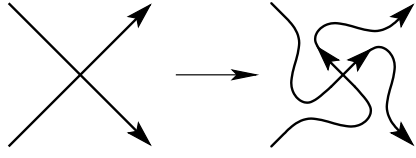}}
\caption{Small perturbation of a divide near a double point}
\label{fig:Xing-turning}
\end{figure}

\begin{proof}
Recall that divides under consideration are additionally equipped with an
orientation. Without a loss of generality, I assume that at every double point
the two branches of $P$ are oriented upwards. One can always perturb $P$
slightly without changing its isotopy type in order to satisfy this condition,
as it is shown in Figure~\ref{fig:Xing-turning}. I also assume that the
inverse self-tangency perestroika is happening ``horizontally'', i.e. the two
new double points are created side by side like depicted in
Figure~\ref{fig:divide-invtangency}. Denote these two points by $a$ and $b$.

\begin{figure}
\centerline{\input{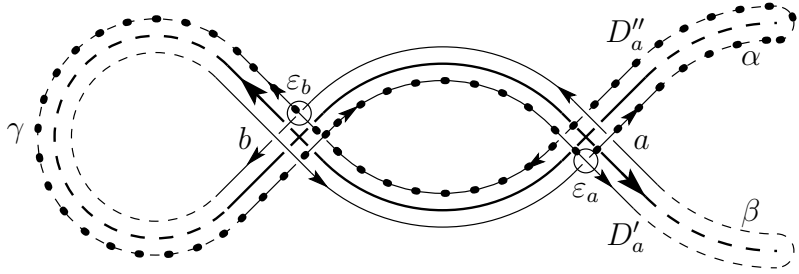}}
\caption{Positive inverse self-tangency perestroika of a divide}
\label{fig:divide-invtangency}
\end{figure}

The common part of $P$ and $P'$ consists of $3$ non-closed curves. Two of them
contain initial and final points of $P$ (or $P'$) and I denote them by $\Ga$
and $\Gb$, respectively. The third curve is denoted by $\Gg$. Without a loss of
generality, one may assume that $\Ga$, $\Gb$, and $\Gg$ are placed as in
Figure~\ref{fig:divide-invtangency}. They may, of course, intersect each other
and have multiple self-intersections.

Observe now that there are two crossing points of $L(P')$ close to $a$ and $b$
such that changing over- and under-crossings at them allows one to pull $4$
strands of $L(P')$ placed between $a$ and $b$ away from each other and to
obtain $L(P)$. Those points are marked with small circles in
Figure~\ref{fig:divide-invtangency} and are denoted by $\Ge_a$ and $\Ge_b$,
respectively. More precisely, let $K_a$ be the result of over- and
under-crossings change of $L(P')$ at $\Ge_a$ and let $K_b$ be the result of
such a change of $K_a$ at $\Ge_b$. Then $K_b$ and $L(P)$ are ambiently
isotopic and, hence, $v_2(L(P))=v_2(K_b)$. It is well known how the value of
$v_2$ behaves under these crossing changes. Namely, the following skein
relation holds true:
\begin{equation}\label{eq:skein}
v_2\left(\vcenter{\hbox{\noindent\begin{picture}(0,0)%
\includegraphics{pos_Xing.pstex}%
\end{picture}%
\setlength{\unitlength}{3947sp}%
\begingroup\makeatletter\ifx\SetFigFont\undefined
\def\x#1#2#3#4#5#6#7\relax{\def\x{#1#2#3#4#5#6}}%
\expandafter\x\fmtname xxxxxx\relax \def\y{splain}%
\ifx\x\y   
\gdef\SetFigFont#1#2#3{%
  \ifnum #1<17\tiny\else \ifnum #1<20\small\else
  \ifnum #1<24\normalsize\else \ifnum #1<29\large\else
  \ifnum #1<34\Large\else \ifnum #1<41\LARGE\else
     \huge\fi\fi\fi\fi\fi\fi
  \csname #3\endcsname}%
\else
\gdef\SetFigFont#1#2#3{\begingroup
  \count@#1\relax \ifnum 25<\count@\count@25\fi
  \def\x{\endgroup\@setsize\SetFigFont{#2pt}}%
  \expandafter\x
    \csname \romannumeral\the\count@ pt\expandafter\endcsname
    \csname @\romannumeral\the\count@ pt\endcsname
  \csname #3\endcsname}%
\fi
\fi\endgroup
\begin{picture}(324,324)(589,-73)
\end{picture}
}}\!\!\right)-
v_2\left(\vcenter{\hbox{\noindent\begin{picture}(0,0)%
\includegraphics{neg_Xing.pstex}%
\end{picture}%
\setlength{\unitlength}{3947sp}%
\begingroup\makeatletter\ifx\SetFigFont\undefined
\def\x#1#2#3#4#5#6#7\relax{\def\x{#1#2#3#4#5#6}}%
\expandafter\x\fmtname xxxxxx\relax \def\y{splain}%
\ifx\x\y   
\gdef\SetFigFont#1#2#3{%
  \ifnum #1<17\tiny\else \ifnum #1<20\small\else
  \ifnum #1<24\normalsize\else \ifnum #1<29\large\else
  \ifnum #1<34\Large\else \ifnum #1<41\LARGE\else
     \huge\fi\fi\fi\fi\fi\fi
  \csname #3\endcsname}%
\else
\gdef\SetFigFont#1#2#3{\begingroup
  \count@#1\relax \ifnum 25<\count@\count@25\fi
  \def\x{\endgroup\@setsize\SetFigFont{#2pt}}%
  \expandafter\x
    \csname \romannumeral\the\count@ pt\expandafter\endcsname
    \csname @\romannumeral\the\count@ pt\endcsname
  \csname #3\endcsname}%
\fi
\fi\endgroup
\begin{picture}(324,324)(589,-73)
\end{picture}
}}\!\!\right)=
lk\left(\vcenter{\hbox{\noindent\begin{picture}(0,0)%
\includegraphics{smooth_Xing.pstex}%
\end{picture}%
\setlength{\unitlength}{3947sp}%
\begingroup\makeatletter\ifx\SetFigFont\undefined
\def\x#1#2#3#4#5#6#7\relax{\def\x{#1#2#3#4#5#6}}%
\expandafter\x\fmtname xxxxxx\relax \def\y{splain}%
\ifx\x\y   
\gdef\SetFigFont#1#2#3{%
  \ifnum #1<17\tiny\else \ifnum #1<20\small\else
  \ifnum #1<24\normalsize\else \ifnum #1<29\large\else
  \ifnum #1<34\Large\else \ifnum #1<41\LARGE\else
     \huge\fi\fi\fi\fi\fi\fi
  \csname #3\endcsname}%
\else
\gdef\SetFigFont#1#2#3{\begingroup
  \count@#1\relax \ifnum 25<\count@\count@25\fi
  \def\x{\endgroup\@setsize\SetFigFont{#2pt}}%
  \expandafter\x
    \csname \romannumeral\the\count@ pt\expandafter\endcsname
    \csname @\romannumeral\the\count@ pt\endcsname
  \csname #3\endcsname}%
\fi
\fi\endgroup
\begin{picture}(324,324)(589,-73)
\end{picture}
}}\!\!\right),
\end{equation}
where $lk$ is the linking number of a two-component link.

Let now $L_a$ be the result of smoothing of $L(P')$ at $\Ge_a$ with respect to
the orientation and let $L_b$ be the result of smoothing of $K_a$ at $\Ge_b$.
Both $L_a$ and $L_b$ are links with two components $L'_a$, $L''_a$ and $L'_b$,
$L''_b$, respectively. Then \eqref{eq:skein} implies that
$v_2(L(P'))-v_2(K_a)=lk(L'_a,L''_a)$ and $v_2(K_a)-v_2(K_b)=lk(L'_b,L''_b)$.

Thus one concludes that $\GD v_2= lk(L'_a,L''_a) + lk(L'_b,L''_b)$. Let us
compute the first linking number. Denote the diagrams of $L'_a$ and $L''_a$ by
$D'_a$ and $D''_a$, respectively. Then the linking number in question
is half the sum of signs of crossings from $D'_a\cap D''_a$. Here the {\em
sign of a crossing} is $+1$ or $-1$ if it looks like the first or second
summand in~\eqref{eq:skein}, respectively. In
Figure~\ref{fig:divide-invtangency} $D''_a$ is depicted with a dotted line.

\begin{figure}
\centerline{\input{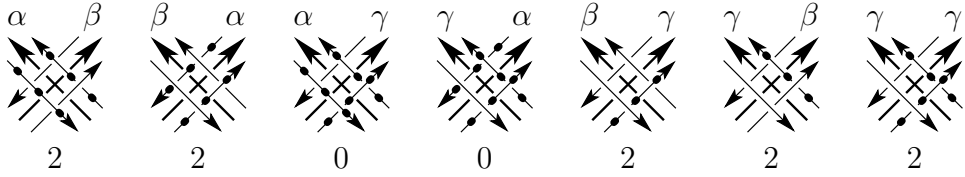}}
\caption{Contribution of crossing points to $2lk(L'_a,L''_a)$}
\label{fig:crossing-contributions}
\end{figure}

\begin{figure}
\centerline{\input{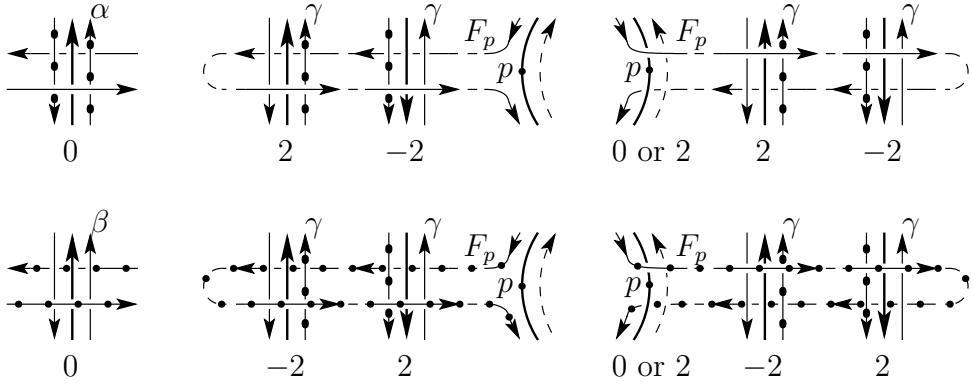}}
\caption{Contribution of ``jumps through infinity'' to $2lk(L'_a,L''_a)$}
\label{fig:minmax-contributions}
\end{figure}

$D'_a$ and $D''_a$ may intersect each other either near a double point of the
divide or where a branch of $D(P')$ jumps through infinity.
Self-intersections of $\Ga$ and $\Gb$ do not contribute to twice the linking
number, since the corresponding pieces of $D(P)$ belong entirely to $D''_a$
and $D'_a$, respectively. Contribution of other double points is shown in
Figure~\ref{fig:crossing-contributions}. It is $2$ for every intersection
point of $\Gb$ with $\Ga$ and $\Gg$ and every self-intersection of $\Gg$. The
contribution is $0$ for every intersection of $\Ga$ with $\Gg$. Points $a$ and
$b$ both contribute $2$. Hence the part of $lk(L'_a,L''_a)$ coming from
crossings located near self-intersections of $P'$ is $\#(\Ga\cap\Gb) +
\#(\Gb\cap\Gg) + n(\Gg) + 2$.

Let $p$ be a point of local minimum or maximum on $P'$ with respect to the
projection onto the $x$-axis. Denote by $F_p$ the corresponding piece of
$D(P')$ that makes a jump through infinity. One may assume that all these
pieces are pairwise disjoint and intersect $P'$ transversely at finitely many
points. Intersections of $F_p$ with branches of $D(P')$ corresponding to
$\Ga$ and $\Gb$ do not contribute to twice the linking number (see the two
leftmost pictures in Figure~\ref{fig:minmax-contributions}). Contribution from
intersections with $\Gg$ is either $2$ if $F_p$ belongs to $D'_a$ ($D''_a$)
and the intersected branch of $\Gg$ is oriented upwards (downwards) or $-2$
otherwise (see Figure~\ref{fig:minmax-contributions}). Moreover, one should
add $2$ if $p\in\Gg$ and is a local maximum. All in all, the contribution of
$F_p$ to twice the linking number is summarized in the following table

\begin{center}
\def\arraystretch{1.3}
\def\phm{\phantom{-}}
\begin{tabular}{|c|c||c|c|}
\hline
Location of $p$&Contribution of $F_p$&
Location of $p$&Contribution of $F_p$\\
\hline
$p\in\minimax{\Ga}mS$&$\phm2\ind_\barg(p)$&
$p\in\minimax{\Gg}mU$&$   -2\ind_\barg(p)-1$\\
$p\in\minimax{\Ga}MS$&$   -2\ind_\barg(p)$&
$p\in\minimax{\Gg}mD$&$\phm2\ind_\barg(p)-1$\\
$p\in\minimax{\Gb}mS$&$   -2\ind_\barg(p)$&
$p\in\minimax{\Gg}MU$&$\phm2\ind_\barg(p)+1$\\
$p\in\minimax{\Gb}MS$&$\phm2\ind_\barg(p)$&
$p\in\minimax{\Gg}MD$&$   -2\ind_\barg(p)+1$\\
\hline
\end{tabular}
\end{center}

\noindent
where $\barg$ is the natural closure of $\Gg$ and $\minimax{C}mS$ and
$\minimax{C}MS$ denote the sets of all local minima and maxima on the curve
$C$, respectively (with an appropriate orientation, if necessary). Recall that
$\ind_\Gg(p)$ is half-integer for $p\in\Gg$. Hence the contribution of $F_p$
is always even, as it should be.

It is obvious that
$\#\big(\,\minimax{\barg}MS\big)-\#\big(\minimax{\barg}mS\,\big)=0$, since
$\barg$ is closed. The difference between $\Gg$ and $\barg$ is a local
maximum near the point $b$. Therefore
$\#\big(\,\minimax{\Gg}MS\big)-\#\big(\minimax{\Gg}mS\,\big)=-1$. Finally

\begin{equation*}
\begin{split}
lk(L'_a,L''_a)=
&\sum_{p\in\minimax{\Ga}mS\cup\minimax{\Gb}MS}\ind_\barg(p)
-\sum_{p\in\minimax{\Gb}mS\cup\minimax{\Ga}MS}\ind_\barg(p)
+\sum_{p\in\minimax{\Gg}mD\cup\minimax{\Gg}MU}\ind_\barg(p)
-\sum_{p\in\minimax{\Gg}mU\cup\minimax{\Gg}MD}\ind_\barg(p)\\
&+\#(\Ga\cap\Gb)+\#(\Gb\cap\Gg)+n(\Gg)+3/2.
\end{split}
\end{equation*}

Computation of $lk(L'_b,L''_b)$ is almost the same. One should only exchange
the roles of $\Ga$ and $\Gb$ everywhere and to take into account that $a$
does not contribute to the linking number anymore. The contribution of $b$
remains $2$. Therefore

\begin{equation}\label{eq:Dv2-inverse}
\begin{split}
\GD v_2=
&\;2\!\!\!\!\sum_{p\in\minimax{\Gg}mD\cup\minimax{\Gg}MU}\ind_\barg(p)
-2\!\!\!\!\sum_{p\in\minimax{\Gg}mU\cup\minimax{\Gg}MD}\ind_\barg(p)\\
&+2\#(\Ga\cap\Gb)+\#(\Ga\cap\Gg)+\#(\Gb\cap\Gg)+2n(\Gg)+2.
\end{split}
\end{equation}

It is simpler to compute $\GD\ScrX$. First of all, $J^+(\overline
P)=J^+(\overline {P'})$ and $\St(\overline P)=\St(\overline {P'})$, since
neither $J^+$ nor $\St$ changes under an inverse self-tangency perestroika.
Let $v$ be a double point of $P$. If $v$ is a self-intersection point of either
$\Ga$, $\Gb$ or $\Gg$, then neither $\tJ(O_v)$ nor $\#(O_v\cap I_v)$ changes.
If $v\in\Ga\cap\Gg$ or $v\in\Gb\cap\Gg$, then $\tJ(O_v)$ does not change, but
$\#(O_v\cap I_v)$ increases by $2$ under the perestroika. Finally, if
$v\in\Ga\cap\Gb$, then $\#(O_v\cap I_v)$ does not change, but $O_v$
experiences a positive inverse self-tangency perestroika and, hence,
$\tJ(O_v)$ increases by $2$. Moreover, $\#(O_a\cap I_a)=\#(O_b\cap
I_b)=\#(\Ga\cap\Gg)+\#(\Gb\cap\Gg)$. Combining these facts together one can
conclude that

\begin{equation}\label{eq:DX-inverse}
\GD\ScrX=2\#(\Ga\cap\Gb)+\#(\Ga\cap\Gg)+\#(\Gb\cap\Gg)+\tJ(O_a)+\tJ(O_b).
\end{equation}

Subtracting \eqref{eq:DX-inverse} from \eqref{eq:Dv2-inverse} one gets

\begin{equation*}
\GD v_2-\GD\ScrX=
2\!\!\!\!\sum_{p\in\minimax{\Gg}mD\cup\minimax{\Gg}MU}\ind_\barg(p)
-2\!\!\!\!\sum_{p\in\minimax{\Gg}mU\cup\minimax{\Gg}MD}\ind_\barg(p)+2n(\Gg)+2
-\tJ(O_a)-\tJ(O_b).
\end{equation*}

It easily follows from Lemma~\ref{lem:J-minmax-formula} below that this
expression is always zero. Indeed, $\tJ(O_a)=\tJ(\Gg)+1-(\ind_\barg(a)-1/2)$
and $\tJ(O_b)=\tJ(\Gg)+\ind_\barg(b)+1/2$ with $\ind_\barg(a)=\ind_\barg(b)$,
where $\tJ(\Gg)$ denotes
$\displaystyle\sum_{p\in\minimax{\Gg}mD\cup\minimax{\Gg}MU}\ind_\barg(p)
-\sum_{p\in\minimax{\Gg}mU\cup\minimax{\Gg}MD}\ind_\barg(p)+n(\Gg)$.

\end{proof}

\begin{lem}\label{lem:J-minmax-formula}
Let $C$ be an oriented closed plane curve such that
\begin{itemize}
\item at every double point the two branches of $C$ are parallel to the main
diagonals $y=\pm x$ of the plane;
\item $C$ has only finitely many points where the tangent vector is parallel
to the $y$-axis, and at every such a point the projection of $C$ onto the
$x$-axis has either a local minimum or a local maximum.
\end{itemize}
Then $\displaystyle\tJ(C)=\#(\vertexup{C})+\#(\vertexdown{C})+
\sum_{p\in\minimax{C}mD\cup\minimax{C}MU}\ind_C(p)
-\sum_{p\in\minimax{C}mU\cup\minimax{C}MD}\ind_C(p).$
\end{lem}
The proof is elementary and is left to the reader. Hint: smoothing at a double
point of $C$ with respect to the orientation does not change the expression
above.

\begin{figure}
\centerline{\input{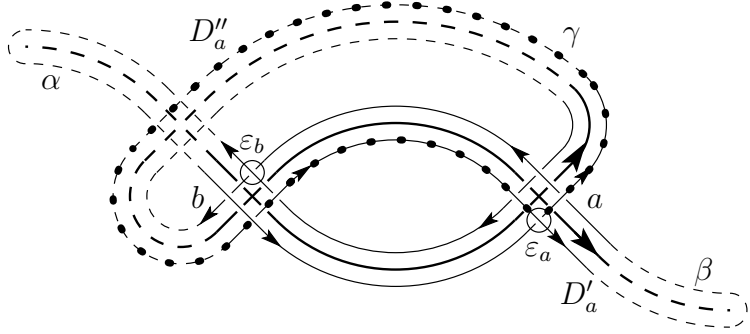}}
\caption{Positive direct self-tangency perestroika of a divide}
\label{fig:divide-dirtangency}
\end{figure}

\begin{lem}\label{lem:direct-move}
Let $P'$ be the result of a positive direct self-tangency perestroika of $P$.
Assign $\GD\ScrX=\ScrX(P')-\ScrX(P)$ and $\GD v_2=v_2(L(P'))-v_2(L(P))$. Then
$\GD\ScrX=\GD v_2$.
\end{lem}
The proof is similar to the one of Lemma~\ref{lem:inverse-move}. The
corresponding picture of a divide and its link close to the place of
perestroika is shown in Figure~\ref{fig:divide-dirtangency}. In this case
\begin{equation}\label{eq:DX-direct}
\begin{split}
\GD v_2
&=2\!\!\!\!\sum_{p\in\minimax{\Gg}mD\cup\minimax{\Gg}MU}\ind_\barg(p)
-2\!\!\!\!\sum_{p\in\minimax{\Gg}mU\cup\minimax{\Gg}MD}\ind_\barg(p)
+\#(\Ga\cap\Gg)+\#(\Gb\cap\Gg)+2n(\Gg)+1\\
&=\#(\Ga\cap\Gg)+\#(\Gb\cap\Gg)+\tJ(O_a)+\tJ(O_b)+1=\GD\ScrX.
\end{split}
\end{equation}

Turn attention to triple point perestroikas now. It follows from the
definition that if a divide $P$ experiences a triple point perestroika, then
the ambient isotopy type of the corresponding link does not change.

\begin{figure}
\centerline{\input{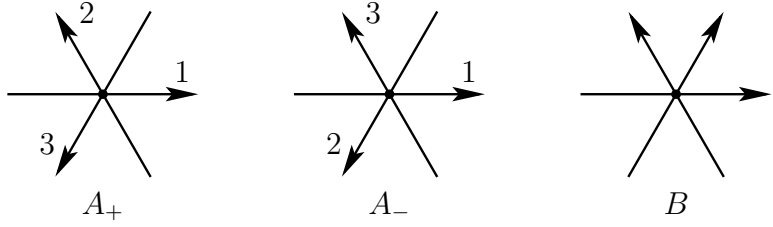}}
\caption{Different types of triple points}
\label{fig:triple-types}
\end{figure}

\begin{figure}
\centerline{\input{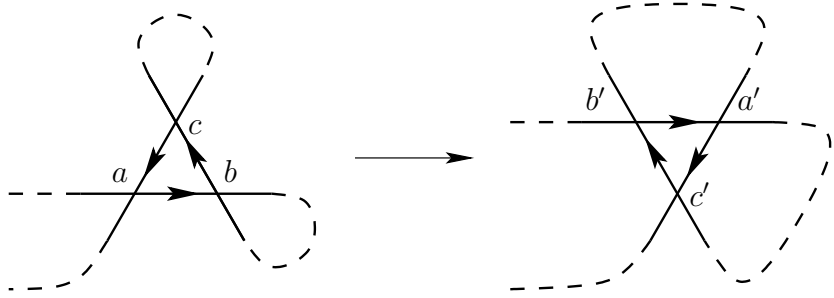}}
\caption{Positive triple point perestroika of type $A_+$}
\label{fig:divide-triple}
\end{figure}

\begin{lem}\label{lem:triple-move}
Let $P'$ be the result of a positive triple point perestroika of $P$. Assign
$\GD\ScrX=\ScrX(P')-\ScrX(P)$. Then $\GD\ScrX=0$.
\end{lem}
\begin{proof}
One can distinguish two kinds of triple points, depending on whether the three
branches of a curve (or a divide) at the point are directed into the same
half-plane or not. A triple point is said to be of {\em type $B$} in the
former case and of {\em type $A$} in the latter one (see
Figure~\ref{fig:triple-types}). Accordingly, there are two kinds of triple
point perestroikas. It is easy to check that a triple point perestroika of
type $B$ is equivalent to a sequence of a type $A$ triple point perestroika
and several self-tangency perestroikas. It is therefore enough to consider
only type $A$ perestroikas in the proof.

Triple points of type $A$ can further be classified as being of type $A_+$ or
$A_-$, depending of whether the cyclic order of the branches defines a
positive or negative orientation of the plane (see
Figure~\ref{fig:triple-types}). Without a loss of generality, I will restrict
my consideration to perestroikas of type $A_+$ only. The corresponding picture
of divides before and after a perestroika is shown in
Figure~\ref{fig:divide-triple}.

It is now straightforward to see that $\St(\overline {P'})=\St(\overline
P)+1$, $\tJ(O_{a'})=\tJ(O_a)-2$, $\tJ(O_{b'})=\tJ(O_b)$,
$\tJ(O_{c'})=\tJ(O_c)$, and $\#(O_{v'}\cap I_{v'})=\#(O_v\cap I_v)+2$ for
$v\in\{a, b, c\}$. Hence $\GD\ScrX=1/2-2+6/4=0$.
\end{proof}

\subsection{Applications to Chmutov's $J^\pm_2$ invariant}\label{sec:chmutov}
$J^\pm_2$ is a second order $J^\pm$-type invariant of long curves (or
I-divides). Chmutov~\cite{Chmutov-J_pm} defined it by explicitly specifying
its actuality table (i.e. values of the invariant on all the chord diagrams
with two chords) and its values on standard divides with at most one
self-tangency point. He also proved that $v_2(L(P))=J^\pm_2(P)$ for any
I-divide $P$. The definition of $J^\pm_2$ did not allow one to integrate this
invariant, i.e. to compute its values on chord diagrams with one chord. Since
these values are nothing more than changes of $J^\pm_2$ under self-tangency
perestroikas, formulas~\eqref{eq:DX-inverse} and~\eqref{eq:DX-direct} provide
an answer to this question.

\begin{cor}
Let $P$ be an I-divide and $P'$ be the result of a positive inverse
self-tangency perestroika of $P$. Let $\Ga$, $\Gb$, $a$, and $b$ be as
in the proof of Lemma~\ref{lem:inverse-move} (see
Figure~\ref{fig:divide-invtangency}). Then
\begin{equation}
\begin{split}
J^\pm_2(P')-J^\pm_2(P)
&=2\tJ(O_a)+2\ind_{O_a}(b)+\#(O_a\cap I_a))+2\#(\Ga\cap\Gb)-1\cr
&=2\tJ(O_b)-2\ind_{O_b}(a)+\#(O_b\cap I_b))+2\#(\Ga\cap\Gb)+1.\cr
\end{split}
\end{equation}
\end{cor}

The proof follows from~\eqref{eq:DX-inverse} and the facts that $\#(O_a\cap
I_a)=\#(O_b\cap I_b)=\#(\Ga\cap\Gg)+\#(\Gb\cap\Gg)$,
$\tJ(O_a)=\tJ(O_b)-2\ind_{O_a}(b)+1$ and $\ind_{O_a}(b)=\ind_{O_b}(a)$.

\begin{cor}
Let $P$ be an I-divide and $P'$ be the result of a positive direct
self-tangency perestroika of $P$. Let $a$ and $b$ be as in the proof of
Lemma~\ref{lem:direct-move} (see Figure~\ref{fig:divide-dirtangency}). Then
\begin{equation}
J^\pm_2(P')-J^\pm_2(P)=2\tJ(O_a)+\#(O_a\cap I_a))=2\tJ(O_b)+\#(O_b\cap I_b)).
\end{equation}
\end{cor}

The proof follows from~\eqref{eq:DX-direct} and the facts that $\#(O_a\cap
I_a)=\#(O_b\cap I_b)=\#(\Ga\cap\Gg)+\#(\Gb\cap\Gg)+1$ and $\tJ(O_a)=\tJ(O_b)$.


\bigskip

\begin{thebibliography}{99999}
\bibitem{ACampo-divide-def}
N. A'Campo, {\sl Le Groupe de Monodromie du D\'eploiement des
    Singularit\'es Isol\'ees de Courbes Planes II}, Actes du Con\-gr\`es
    Inter\-national des Math\'ema\-ti\-ciens, tome 1, 395--404, Vancouver,
    1974.

\bibitem{ACampo-divide-1}
N. A'Campo, {\sl Real deformations and complex topology of plane curve
    singularities}, Ann. Fac.  Sc. de Toulouse Math. (6){\bf 8} (1999), 5--23,
    arXiv: alg-geom/9710023.

\bibitem{ACampo-divide-2}
N. A'Campo, {\sl Generic immersions of curves, knots, monodromy and Gordian
    number}, Publ. Math. Inst. Hautes \'Etudes Sci. {\bf 88} (1998),
    151--169, (1999), arXiv: math.GT/9803081.

\bibitem{ACampo-divide-3}
N. A'Campo, {\sl Planar trees, slalom curves and hyperbolic knots}, Publ.
    Math. Inst. Hautes \'Etudes Sci. {\bf 88} (1998), 171--180, (1999),
    arXiv: math.GT/9906087.

\bibitem{Aicardi-tree-curves}
F.~Aicardi, {\sl Tree-like curves}, In: {\sl Singularities and Bifurcations},
    Adv. Sov. Math. {\bf 21} (1994), AMS, Providence, RI, 1--31.

\bibitem{Arnold-curve-invars}
V. I. Arnold, {\sl Plane curves, their invariants, perestroikas
    and classifications}, In: {\sl Singularities and Bifurcations},
    Adv. Sov. Math. {\bf 21} (1994), AMS, Providence, RI, 33--91.

\bibitem{Chmutov-J_pm}
S. Chmutov, {\sl An integral generalization of the Gusein-Zade--Natanzon
    theorem}, arXiv: math.GT/0207238.

\bibitem{Chmuzhin-invars}
S. Chmutov, S. Duzhin, {\sl Explicit formulas for Arnold's generic
    curve invariants}, In: {\sl The Arnold-Gelfand mathematical seminars},
    Birkh\"auser Boston, Boston, MA, 1997, 123--138.

\bibitem{Gibson-Ishikawa-diagram}
W. Gibson, M. Ishikawa, {\sl Links of oriented divides and fibrations in link
    exteriors}, to appear in Osaka J. Math.

\bibitem{Hirasawa-diagram}
M. Hirasawa, {\sl Visualization of A'Campo's fibered links and unknotting
    operations}, Proceedings of the First Joint Japan-Mexico Meeting in
    Topology (Morelia, 1999). Topology Appl. {\bf 121} (2002), no. 1--2,
    287--304.

\bibitem{Polyak-arnold-gauss}
M. Polyak, {\sl Invariants of curves and fronts via Gauss diagrams}, Topology
    {\bf 37} (1998), no. 5, 989--1009.

\bibitem{Myself-strangeness}
A. Shumakovitch, {\sl Explicit formulas for strangeness of plane
    curves}, Algebra i Analiz {\bf 7} (1995), no. 3, 165--199; English transl.
    in St.~Petersburg Math. J. {\bf 7} (1996), no. 3, 445--472.


\bibitem{Viro-J's-formulas}
O. Ya. Viro, {\sl First degree invariants of generic curves on surfaces},
    Preprint, Dept. Math., Uppsala Univ., 1994:21.

\end{thebibliography}
\end{document}